\documentclass[a4paper,12pt,leqno]{amsart}
\usepackage{graphicx}
\usepackage{enumerate}

\catcode`\@=11
\newcommand{\languefrancaise}{
\catcode `\?=\active
\def?{\relax\ifhmode\ifdim\lastskip>\z@\unskip\fi
\kern.2em\fi
\string?}

\catcode `\;=\active
\def;{\relax\ifhmode\ifdim\lastskip>\z@\unskip\fi
\kern.2em\fi
\string;}

\catcode `\:=\active
\def:{\relax\ifhmode\ifdim\lastskip>\z@\unskip\fi
\penalty\@M\ \fi
\string:}

\catcode `\!=\active
\def!{\relax\ifhmode\ifdim\lastskip>\z@\unskip\fi
\kern.2em\fi
\string!}

\frenchspacing
}
\catcode`\@=12

\long\def\xcom#1{}

\newcommand{\ladate}{\space\the\day\ \ifcase\month\or janvier\or f\'evrier
\or mars\or avril\or mai\or juin\or juillet\or aout\or septembre
\or octobre\or novembre\or d\'ecembre\fi
\ {\oldstyle\the\year}}

\def\alphabet#1{\ifcase#1\or a\or b\or c\or d\or e\or f\or g\or 
h\or i\or j\or k\or l\or m\or n\or o\or p\or q\or r\or s\or t\or 
u\or v\or w\or x\or y\or z\fi}

\newcommand{\up}[1]{\raise 1ex\hbox{\sevenrm#1}}

\newcommand{\R}{\mathbb{R}}
\newcommand{\PP}{\mathbb{P}}
\newcommand{\E}{\mathbb{E}}
\newcommand{\N}{\mathbb{N}}

\newcommand{\cov}{\mathop{\rm Cov}}

\newcommand{\unsur}[1]{{\frac{1}{#1}}}
\def\un#1{{1_{({#1})}}}
\def\valabs#1{{\left\vert {#1} \right \vert}}
\def\esp#1{{{\E}\left [ {#1} \right ]}}
\def\etp#1{{\left ( {#1} \right )}}
\def\etc#1{{\left [ {#1} \right ]}}
\def\crochet#1{{\left < {#1} \right >}}

\def\prob#1{{{\PP}\left ( {#1} \right )}}
\def\norme#1{{\left \Vert #1 \right \Vert}}
\def\ens#1{{\left\{#1\right\}}}

\newcommand{\egaloi}{{\,\mathrel{\mathop{=}\limits^d}\,}}

\def\dessus#1#2{\mathord{\mathop{\kern 0pt #2}\limits^#1}}

\def\convergeenloi#1{\mathrel{
                 \mathop{\longmapsto}\limits^d_{#1\to\infty} 
}}

\newcommand{\cvloit}{\convergeenloi{t}}

\newcommand{\cvt}{\mathrel{\mathop{\longmapsto}\limits_{t\to\infty}}}
\newcommand{\cvn}{\mathrel{\mathop{\longmapsto}\limits_{n\to\infty}}}

\newcommand{\bit}{\begin{itemize}}
\newcommand{\eit}{\end{itemize}}
\newcommand{\ben}{\begin{enumerate}}
\newcommand{\een}{\end{enumerate}}

\newcounter{moncompteur}

\newenvironment{myenumerate}%
{\begin{list}{\arabic{moncompteur}. }{\usecounter{moncompteur}%
\setlength{\leftmargin}{0pt}%
\setlength{\labelwidth}{0pt}%
\setlength{\listparindent}{0pt}%
\setlength{\labelsep}{0pt}}}%
{\end{list}}

\def\bmen{\begin{myenumerate}}
\def\emen{\end{myenumerate}}

\newcommand{\intof}{{\int_0^\infty}}
\newcommand{\tgo}{{\, ; t \geq 0\,}}
\newcommand{\intot}{{\int_0^t}}
\newcommand{\rp}{{\R_+}}

\newcommand{\undemi}{\frac{1}{2}}

\def\a{\alpha}

\newcommand{\Arond}{{\mathcal A}}

\newcommand{\Frond}{{\mathcal F}}

\newcommand{\Irond}{{\mathcal I}}

\newcommand{\Nrond}{{\mathcal N}}

\newcommand{\Yrond}{{\mathcal Y}}

\RequirePackage{xspace}
\newcommand{\ou}{Ornstein-Uhlenbeck\xspace}

\newcommand{\fbm}{fractional Brownian motion\xspace}

\newcommand{\sde}{stochastic differential equation\xspace}

\newcommand{\eps}{\epsilon}

\newtheorem{theorem}{Theorem}
\newtheorem{lemma}[theorem]{Lemma}
\newtheorem{proposition}[theorem]{Proposition}

\theoremstyle{definition}

\theoremstyle{remark}

\begin{document}
\title{Fractional Brownian Motion and The Markov Property}

\author[P. Carmona]{Philippe Carmona}
\address{Philippe Carmona\\
Laboratoire de Statistique et Probabilit\'es\\
Universit\'e Paul Sabatier\\
118 route de Narbonne\\ 31062 Toulouse Cedex 4}
\email{carmona@cict.fr}
\author[L. Coutin]{Laure Coutin}
\address{Laure Coutin\\
Laboratoire de Statistique et Probabilit\'es\\
Universit\'e Paul Sabatier\\
118 route de Narbonne\\ 31062 Toulouse Cedex 4}
\email{coutin@cict.fr}



\newcommand{\bm}{Brownian motion\xspace}
\newcommand{\as}{almost surely\xspace}
\newcommand{\rv}{random variable\xspace}

\begin{abstract}
Fractional \bm belongs to a class
of  long memory Gaussian processes that can be represented 
 as linear functionals of an infinite dimensional Markov process. This  leads naturaly to~:
\begin{itemize}
\item An efficient algorithm to approximate the process.
\item An ergodic theorem  which applies to functionals of the type 
$$\intot \phi(V_h(s))\,ds \quad\text{where}\quad V_h(s)=\intot h(t-u)\, dB_u\,.$$
\end{itemize}

\end{abstract}

\date{\today}\vfuzz=4pt \hfuzz=40pt
\maketitle

\section{Introduction}

Fractional \bm of Hurst parameter $H\in(0,1)$ is a real centered Gaussian process $(W^H_t\tgo)$ of covariance
\begin{equation}
\esp{W^H_t W^H_s} = \frac{c(H)}{2} (s^{2H} + t^{ 2H} - \valabs{t-s}^{2H})\qquad(s,t\ge 0)\,.
\end{equation}
One easily sees that it starts from zero, $W^H_0=0$ \as, it has stationary increments $(W^H_{t+s}-W^H_s\tgo)\egaloi (W^H_t\tgo)$ since

\begin{equation}
\esp{(W^H_t -W^H_s)^2}= c(H) \valabs{t-s}^{2H}\qquad(s,t\ge 0)\,.
\end{equation}

Fractional \bm enjoys the following scaling (self-similarity) property

\begin{equation}
(\unsur{a^H} W^H_{at}\tgo)\egaloi (W^H_t\tgo)\qquad(a>0)\,.
\end{equation}

The increments are independent if and only if $H=\undemi$ (\bm). If $H<\undemi$, they are negatively correlated, and if $H>\undemi$ they are positively correlated. Furthermore it is a long memory process, since, for $H\neq \undemi$,
 the covariance between far apart increments decrease to zero as a power law:
$$ \cov(W^H(t)-W^H(0), W^H(t+u)-W^H(u)) \sim_{u\to +\infty} constant \times u^{2H-2}\,.$$

Fractional \bm has been used as a model of signal or noise in various domains: geophysical data~\cite{graf83,beranterrin96}, communication processes~\cite{mandelbrot65,lelandtww94}, see also the references therein. Therefore it is important to provide good numerical approximating schemes.

\subsection*{Problem I: Approximation of \fbm}

We shall not linger on the first (naive) idea that may cross our mind. We simulate  an $n$-dimensional Gaussian random vector with the same disribution as
$(W^H_{s_1}, \ldots , W^H_{s_n})$ and we perform a linear or spline interpolation. This amounts to computing the square root of the $n\times n$ covariance matrix and is computationnaly inneficient.

A second solution is to use the representation introduced by Mandelbrot and Van Ness~\cite{mandelbrotvanness68}: if $\a=H+\undemi$,

$$ W^H_t = \unsur{\Gamma(\a)} \int_{-\infty}^t ((t-u)^{\a-1} - (-u)_+^{\a-1})\, dB_u$$
where $B$ is a standard \bm on the real line. Without any loss in generality, we can say that the problem is to approximate the finite, but long, memory part:

$$ V(t)= V_h(t) = \intot h(t-s)\, dB_s$$

where $(B_t\tgo)$ is a standard \bm. The na{\"\i}ve way to do this is to fix a small step $\valabs{\Delta_n}=\unsur{n}$ for the partitions we use
$$ \Delta_n= 0 < s_1=\unsur{n} < s_2=\frac{2}{n} < \cdots \,.$$
Hence, if $t=r/n$, then $V(t)$ is approximated by
\begin{align*}
 V(t) &\simeq \sum h(t-s_i)\, (B_{s_{i+1}} -B_{s_i}) \\
&= \unsur{\sqrt{n}} \sum h(t-s_i) \, \xi_i \\
&= \unsur{\sqrt{n}} (h(\unsur{n}) \xi_r + h(2/n) \xi_{r-1} + \cdots )
\end{align*}
where $\xi_i=\sqrt{n}(B_{s_{i+1}}-B_{s_i})$ are independent identically distributed standard normal random variables.

For the next step, $t+\Delta t= t +\unsur{n}=(r+1)/n$, we have

$$ V(t+\Delta t)\simeq \unsur{\sqrt{n}} ( h(1/n) \xi_{r+1} + h(2/n) \xi_r + \cdots )$$

Therefore, all we have to do is to compute the $h(i/n)$, generate the $\xi_i$, and keep them all in memory for the next step.

Unfortunately, if we want to keep the approximation error low (that is $n$ big), this amounts to a lot of data to keep in memory.

The natural idea is to impose a threshold, that is to compute $h(i/n)$ and generate $\xi$ only for $1\le i\le M_n$. This may work well for smooth functions $h$, but not for
$$ h(u) = \unsur{\Gamma(\a)} u^{\a-1}\un{u>0} $$
where $\undemi<\a<1$, since $h$ explodes near zero, and thus our partition is not dense enough near $t$~!

\xcom{
\begin{center}
\begin{figure}
\includegraphics[width=13cm]{graphedeh.eps}
\caption{The function $s\to h(t-s)$ when $h(u)= \unsur{\Gamma(\a)} u^{\a-1}\un{u>0} $}
\end{figure}

\end{center}
}
Unfortunately this example is important since it is in essence fractional \bm of index $H=\a-\undemi$.

In~\cite{carcoumon1997},Carmona, Coutin and Montseny derived an approximation scheme from a representation formula inspired by the diffusive input-output model of fractional integration developped by Audounet, Montseny adn Mbodje~\cite{audounetmontseny95,audounetmontsenymbodje93}.

More precisely Fubini's stochastic theorem implies that if $h$ has the ``spectral'' representation

$$ h(u) = \int  e^{ -ux} \mu(dx) \qquad(u>0)$$
where $\mu$ is a positive measure on $(0,+\infty)$, finite on compact sets, and such that $\int (1\wedge x^{-\undemi})\,\mu(dx) <+\infty$,
then
$$ V_h(t) = \intof  Y^x_t\, \mu(dx) $$
where $$ Y^x_t =\intot e^{-x(t-s)}\, dB_s\,.$$
Observe that for fixed $x\ge 0$, $(Y^x_t\tgo)$ is an \ou process of parameter $x$, that is a Gaussian semi-martingale Markov process solution of the \sde
$$ dY^x_t = dB_t - x\, Y^x_t \, dt\,.$$

Therefore, $ V_h=\crochet{\Yrond}_\mu$
is a linear functional of the infinite (here function valued) Markov process
$$ \Yrond=(\Yrond_t\tgo)\,,\qquad \Yrond_t=(Y^x_t, x \ge 0)\,.$$

The idea behind the approximation scheme is to first do a spatial approximation

$$ V_h(t) \simeq V^\pi_h(t) = \sum_\pi c_i\,Y^{\eta_i}_y\,,$$

where $\pi=\ens{\eta_i}_i$ is a chosen partition of $K$ converging to identity, and $c_i$ are coefficients depending only on $h$. This approximation is uniformly good with respect to time, that is $\sup_{t\le T} \valabs{V_h(t)- V^\pi_h(t)}$ goes to zero when $\pi$ grows.

The second step is the time dicretization, which can be performed in the following way~:

\begin{equation}\label{intro:regle}
 Y^{\eta_i}(t+\Delta t)= Y^{\eta_i}(t) + \Delta B_t - \eta_i Y^{\eta_i}(t)\,.
\end{equation}

This second step is a classical one : see e.g. the book of Bouleau and Lepingle~\cite{bouleaulepingle}, chapter V, section B,  pages 271--278, where it is stated that to get a precision of $\epsilon$, one has to consider at most a constant $C$ times $N=\log(\unsur{\epsilon})$ time points. Nevertheless, when one needs to do the whole family of approximations~\eqref{intro:regle}, one needs to control precisely the dependence of the constant $C=C(\eta_i,T)$ with respect to the space and time variables. The classical method, combined with Gronwall's lemma, gives an exponential upper bound for the $L^p$-norm of the approximation error 
$C(\eta,T)\le c_1(T,p) \eta \exp(c_2(T,p) \eta^p)$. Therefore, we have to use
another time discretization scheme for which we have polynomial constants (see e.g. \cite{carmonacontrol98} or \cite{coutinmontceny1997}).

The first   step amounts to approximate a linear functional of an infinite Markov process by a linear functional of a finite dimensional Markov process.  The secret of the success of this method is that for a large class of functions $h$, this dimension can be kept low : for a precison of $\epsilon$, the cardinal of the partition $\pi(\epsilon)$ is roughly
$$\Nrond(\epsilon)=\valabs{\pi(\epsilon)} \sim c \sqrt{\unsur{\epsilon}}\, \log (\unsur\epsilon)\,.$$

We do not say that this approximation method is the best, but we want to stress the fact that the Markov property of the underlying process $\Yrond$ allows an algorithm that is very economic in memory allocation. You only have to store
 $ \Nrond(\epsilon)$ values at each time step, then simulate one Gaussian random variable, and then update these values for the next step (according to rule~\eqref{intro:regle}).

\subsection*{The example of fractional \bm}

We can take advantage not only of the Gaussian nature of the processes $Y^x$, but also of their semimartingale nature (via It\^o's formula). We have been able to give exact $L^2$ and almost sure convergence rates for fractional \bm of index $H=\a -\undemi \in(0,\undemi)$. 

Indeed, to obtain the precision
$$ \sup_{t\le T} \norme{V_h(t)-V^\pi_h(t)}_{L^2} \le \epsilon \,,$$
we may choose the compact
$$ K(\epsilon)\simeq \etc{\epsilon^{\unsur{1-\a}}, (\unsur{\eps})^{\frac{2}{2\a +3}}}$$
and a geometric partition

$$ \pi(\epsilon) =\ens{r^{-m},r^{-m+1}, \ldots , r^n}$$
of ratio $r=r(\epsilon)= 1 + c(\a) \sqrt{\epsilon}$.

Furthermore, if $\epsilon_n$ decreases fast enough to $0$,
\begin{equation}\label{intro:eq:vitesseepsilon}
 \epsilon_n = o\etp{ (\log n)^{- \frac{2\a +3}{8(2-\a)}}}
\end{equation}
then, we have the uniform almost sure convergence

$$ \sup_{t\le T} \valabs{V_h(t) - V^{\pi(\epsilon_n)}_h(t)} \cvn 0\qquad\text{almost surely}\,.$$

Theoretical results (see~\cite{carcoumon1997} or \cite{coutinmontceny1997}) show then that to obtain a precision of $\epsilon=10^{ -3}$, we need to consider a few hundred points. Practical results show that forty points are enough. Furthermore, we have been able to check that (almost surely) the numerical approximations stick exponentially fast to the actual random process (and this happens uniformly on compact sets), when $\epsilon_n$ goes faster to $0$ than indicated in~\eqref{intro:eq:vitesseepsilon}.

\subsection*{Problem II : an ergodic theorem}
It is well known (Birkhoff's ergodic theorem) that for a stationary  process $U$, we have the convergence, almost sure and in $L^1$,
\begin{equation}
\unsur{t} \intot \phi(U(s))\, ds \cvt \esp{\phi(U(0))}\,,
\end{equation}
as soon as $\phi$ is an integrable Borel function: $\esp{\valabs{\phi(U(0))}}< +\infty$.

Of course, we first have to check the ergodicity, that is that the invariant $\sigma$-field is almost surely trivial, which is not so easy to do in general (see e.g. Rozanov~\cite{rozanov66}, page ???).

But what we actually want to do is to answer the following question. Given $(V(t)\tgo)$ a real Gaussian process, a priori non stationary, such that 
$V^r=(V(t+r)\tgo)$ converges in distribution to the stationary process $U$. We want to show that, as soon as $\phi$ is an integrable Borel function: $\esp{\valabs{\phi(U(0))}}< +\infty$, we have the convergence, almost sure and in $L^1$,
\begin{equation} \label{intro:cvv}
\unsur{t} \intot \phi(V(s))\, ds \cvt \esp{\phi(U(0))}\,.
\end{equation}
To be more precise, given $h\in L^2(\rp)$ and $(B_t,t\in \R)$ a standard Brownian motion, we let 
\begin{equation}
V(t)=\intot h(t-s)\, dB_s\,,\quad\text{and}\quad U(t)=\int_{-\infty}^t h(t-s)\, dB_s\,.
\end{equation}
We shall show that~\eqref{intro:cvv} holds if $h$ is the Laplace transform 
\begin{equation}
h(u)=\int e^{-ux}\, \mu(dx)\,,
\end{equation}
of a positive measure $\mu$, finite on compact sets and such that
\begin{equation}
\int (1\wedge \unsur{\sqrt{x}})\, \mu(dx) < +\infty\,.
\end{equation}

\def\yrt{{\Yrond_t}}
\def\gyrt{{\Gamma_\yrt}}
\def\yf{{\Yrond_\infty}}
\def\gyf{{\Gamma_\yf}}

\section{The Markov and ergodic properties of the underlying infinite dimensional process}

Let $B=(B_t\tgo)$ be a standard \bm. Given $x\ge 0$, $Y^x=(Y^x_t\tgo)$ is the \ou process driven by $B$ of parameter $x$ (starting from $0$):
$$ Y^x_t = \intot e^{-x(t-s)}\, dB_s\,.$$
It is also the solution of the \sde
$$ dY^x_t = dB_t -x Y^x_t\, dt\,.$$

For every $t\ge 0$, we consider $\Yrond_t=(Y^x_t, x>0)$; it is a centered Gaussian process of covariance
$$ \gyrt(x,y)=\frac{1-e^{-t(x+y)}}{x+y}\qquad(x,y> 0)\,.$$
The natural candidate for a limiting process is $\yf$ the centered Gaussian process of covariance
$$ \gyf(x,y)=\unsur{x+y} \qquad(x,y> 0)\,.$$

\begin{proposition} 
 The process $\Yrond=(\yrt\tgo)$  is an infinite dimensional Markov process.
\end{proposition}
\begin{proof}
The strong Markov property is quite easy to prove. Indeed, let $T$ be an almost surely finite stopping time (with respect to $(\Frond_t\tgo)$ the completed, hence right-continuous, filtration generated by the \bm $B$). Then, $\bar{B}=(\bar{B}_t=B_{t+T}-B_T\tgo)$ is a standard \bm independent of $(B_u,u\le T)$, and $(Y^x_{t+T},\tgo)$ is the solution of the \sde
$$ dY^x_{t+T} = d\bar{B}_t - x\, Y^x_{t+T}\, dt\,.$$
Therefore,
$$ Y^x_{t+T}=e^{-xt}Y^x_T + \intot e^{-x(t-s)}\, d\bar{B}_s\,.$$
Hence, for any bounded measurable function $\psi$,

$$ \esp{ \psi(Y_{t+T}\tgo) \mid \Frond_T} = \E_{\Yrond_T}\etc{\psi(Y_t\tgo)} \,,$$ with

$$ \E_y\etc{\psi(\Yrond_t\tgo)}=\esp{\psi((Y^{x,y(x)},x\ge 0)\tgo)} \,,$$
and $(Y^{x,y}_t\tgo)$ is the \ou process driven by $B$, of parameter $x$, starting from $y$.
\end{proof}

\begin{proposition}\label{mar:pro:cvloiy}
Suppose that $\mu$ is a sigma-finite measure on $(0,+\infty)$ such that,  for a $p>1$,

$$\int \mu(dx)\,\sup(x^{-p/2},x^{-\undemi}) < +\infty\,.$$
 Then,
\begin{enumerate}
\item For every $y \in L^1(\mu)$, under $\PP_y$, we have the convergence in distribution
$$ \Yrond_t \cvloit \Yrond^\infty\,.$$
\item The only bounded invariant functions are constant.
\end{enumerate}
\end{proposition}

\begin{proof}
(1) We shall first prove the convergence of characteristic functions, since this will characterize the distribution of a possible limit (see~\cite{borkar95}), and then prove the tightness of the family of distributions of $(\Yrond_t,\tgo)$.

The characteristic function is defined by
$$\phi_{\yrt}(\chi)=\E_y\etc{e^{i \crochet{\chi,\yrt}}}\qquad(\chi \in L^\infty)\,.$$
Observe that $\crochet{\chi,\yrt}$ is, under $\PP_y$, a Gaussian \rv of mean
$$ \int \chi(x)y(x) e^{-xt}\, \mu(dx) \cvt 0 \quad\text{by dominated convergence}$$
and of variance
\begin{align*}
  \esp{\etp{\int \chi(x) Y^x_t \,\mu(dx)}^2} &= \int\int \chi(x)\chi(y) \mu(dx)\mu(dy) \esp{Y^x_t Y^y_t} \\
&= \int\int \chi(x)\chi(y) \mu(dx)\mu(dy) \frac{1-e^{-t(x+y)}}{x+y} \\
&\cvt  \int\int \chi(x)\chi(y) \frac{\mu(dx)\mu(dy)}{x+y}
\end{align*}
This last integral is finite, Indeed, since $ab\le \undemi(a^1+b^2)$,we can majorize it by
$$ \ldots \le \norme{\chi}_{L^\infty}^2\int\int \undemi\frac{\mu(dx)\mu(dy)}{\sqrt{x}\sqrt{y}} = \undemi\norme{\chi}_{L^\infty}^2 \etp{\int \frac{\mu(dx)}{\sqrt{x}}}^2 < +\infty\,.$$

Therefore,  $\crochet{\chi,\yrt}$ converges in distribution to a centered Gaussian variable with variance
$$ \int\int \chi(x)\chi(y) \frac{\mu(dx)\mu(dy)}{x+y}\,,$$
and this is easily seen to be the variance of $\crochet{\chi,\yf}$.

We shall now prove the tightness of the family of laws of the processes
$(\Yrond_t, t\ge 1) $. We need to introduce a little topology: we consider that for every $t>0$, $\yrt$ is a random variable taking its values in $L^1(\mu)$ endowed with the topology $\sigma(L^1,L^\infty)$; this is the coarsest topology for which the functions 
$$ y\to \crochet{\chi,y}\quad(\chi \in L^\infty(\mu))$$
are continuous. Indeed,
$$\esp{\int \valabs{Y^x_t}\,\mu(dx)}= c \int \mu(dx) \etp{\frac{1-e^{-2xt}}{2x}}^\undemi < +\infty$$
\xcom{
and when $\int \mu(dx) \unsur{\sqrt{x}} <+\infty$, we can prove directly that
$\Yrond_\infty$ belongs to $L^1(\mu)$:
$$\esp{\int \valabs{Y^x_\infty}\,\mu(dx)}= c \int \mu(dx) \etp{\unsur{2x}}^\undemi < +\infty$$
}

Let $\phi(x)=1\wedge\unsur{\sqrt{x}}$. Given $M>0$, Lemma~\ref{mar:lem:ui} 
implies that, if $1=1/p +1/q$, then
$$ K=\ens{y\in L^1: \int \frac{\valabs{y}^p}{\phi^{p/q}}\,d\mu \le M}$$
is a relatively compact set. 

Let us notice that, by construction, under $\PP_y$ the law of $\Yrond_t$ is the law , under $\PP_0$ of $\Yrond_t +(x\to e^{-xt}y(x))$. Since this last function converges weakly to the null function, we can suppose $y=0$, without any loss in generality. We have the sequence of inequalities, for $t\ge 1$,
\begin{align*}
  \prob{\Yrond_t \notin K} &= \prob{ \int \valabs{Y^x_t}^p \unsur{\phi(x)^{p/q}}\, \mu(dx) >M^2} \\
&\le \unsur{M} \esp{\int \valabs{Y^x_t}^p \unsur{\phi(x)^{p/q}}\, \mu(dx)} \\
&\le \unsur{M} \esp{\int \etp{\frac{1-e^{-2xt}}{2x}}^{p/2} \unsur{\phi(x)^{p/q}}\, \mu(dx)} \\
&\le \unsur{M} \esp{\int \etp{\frac{1}{2x}}^{p/2} \unsur{\phi(x)^{p/q}}\, \mu(dx)}\,.
\end{align*}
This last quantity goes uniformly in $t\ge 1$ to zero as $M\to \infty$ , since the integral is finite.

(2) (i) First of all, we have to prove that if $f$ is a bounded function depending on a finite number of coordinate, then $P_t f(y) \cvt \esp{f(\yf)}$. Indeed, we have, for a bounded measurable function $g:\R^n\to \R$, and ${\chi_i, 1\le i\le n}\subset L^\infty(\mu)$
$$ f(y) = g(\crochet{\chi_i,y}, 1\le i\le n)\,.$$

Then
\begin{align*}
P_tf(y)&= \E_y\etc{g(\crochet{\chi_i,\yrt},1\le i\le n)}  \\
&= \esp{g(\crochet{\chi_i,y(.)e^{-t.}} + \crochet{\chi_i,\yrt},1\le i\le n)}\,.
\end{align*}

But we know that if $X=(X_1, \ldots, X_n)$ is a non degenerate centered Gaussian vector, then we have the absolute continuity relationship
$$ \esp{g(z+X)} = \esp{\Lambda(z)g(X)} \qquad(z \in \R^n)\,,$$
with $\Lambda$  a continuous bounded function. Applying this to $(\crochet{\chi_i,\yrt}, 1\le i\le n)$, for $t>0$, yields

$$ P_t f(y) = \esp{\Lambda(\crochet{\chi_i,y(.)e^{-t.}},1\le i\le n) g(\crochet{\chi_i,\yrt},1\le i\le n)}\,.$$

Now if $y_p \to y$ weakly in $L^1(\mu)$, then, for each $i$, 

$$\crochet{\chi_i,y_p(.)e^{-t.}} \to \crochet{\chi_i,y(.)e^{-t.}}\,$$
and $P_tf(y_p)\to P_tf(y)$. Therefore, if $t>0$, then $P_t f$ is a continuous bounded function. Since, under $\PP_y$, $\yrt \cvloit \Yrond_\infty$, the law of $\yf$ is an invariant measure for the semi-group $(P_t\tgo)$, and we have
$$ \lim_{u\to \infty} P_u f=\lim_{s\to \infty} P_s(P_t f)= \esp{P_t f(\yf)}=\esp{f(\yf)}\,.$$

(ii) Now, Let $f:L^1(\mu)\to \R$ be a bounded invariant function; that is, 
$$ f(y)=P_tf(y)=\E_y\etc{f(\Yrond_t)} \qquad(y\in L^1, t\ge 0)\,.$$
Observe that the sigma field $\Frond=\sigma(\chi, \chi \in L^\infty(\mu))$, is generated by the algebra of sets depending on a finite number of coordinate
$$\Arond=\Arond\ens{\chi, \chi\in L^\infty(\mu)}$$
that is the smallest algebra containings the sets $\ens{y:\crochet{\chi,y}\in B}$ , $\chi \in  L^\infty(\mu)$, $B$ borel set of $\R$. Therefore, the measure $P_{\Yrond_\infty}$, the law of $\Yrond_\infty$, being finite, there exists bounded functions $(f_n,g_n,n\in\N)$ depending only on a finite number of coordinate such that
$$ f_n\le g \le g_n \qquad \esp{\valabs{f_n -g_n}(\Yrond_\infty)} \cvn 0\,.$$
We take the limit as $t\to +\infty$  in
$$ P_t f_n(y) \le P_tf(y)=f(y) \le P_tg_n(y) $$
to obtain
$$ \esp{f_n(\yf)} \le f(y) \le \esp{g_n(\yf)}$$
and then let $n$ go to $\infty$ to conclude that $f(y)=\esp{(\yf)}$ is a constant function.

\end{proof} 

\begin{lemma}\label{mar:lem:ui}

  Let $\mu$ be  a sigma-finite measure, and $\phi$ be a strictly positive $\phi >0$, integrable function. Then, for every $M\ge 0$, $p,q>1$ such taht $1=1/p+1/q$,

$$ K=\ens{f: \int \frac{\valabs{f}^{p}}{\phi^{p/q}}\,d\mu \le M}$$
is a relatively compact set of $L^1(\mu)$ endowed with the weak topology $\sigma(L^1,L^\infty)$.
\end{lemma}

\begin{proof}

  Following Theorem~IV.8.9 of Dunford--Schwarz~\cite{dunsch88} we have to prove that 
  \begin{enumerate}
  \item $K$ is a bounded subset of $L^1$ .
\item If $E_n$ is a sequence of measurable sets that decreases to $\emptyset$, then $\sup_{f\in K} \valabs{\int_{E_n} f\, d\mu}\cvn 0$.
  \end{enumerate}
By H\"older's inequality, for $f$ in $K$
\begin{align*}
  \int 1_{E_n}\valabs{f}\, d\mu&\le \norme{\frac{f}{\phi^{p/q}}}_{L^p(\mu)} \norme{1_{E_n} \phi^{1/q}}_{L^q(\mu)} \\
& \le M^{1/p} \etp{\int 1_{E_n} \phi d\mu}^{1/q}\,.
\end{align*}
Taking $E_0=\Omega$ the whole space, we see that $K$ is bounded in $L^1$. The dominated convergence theorem gives that $\int 1_{E_n} \phi\, d\mu \cvn 0$, so the convergence is uniform for $f$ in $K$.
\end{proof}

\xcom{
Let $\gamma_1$ denote the standard Gaussian measure on the line:
$$ \gamma_1(dx)=\unsur{\sqrt{2\pi}}\,e^{-x^2/2}\, dx$$
and $N\sim \Nrond(0,1)$ be a standard Gaussian variable (of law $\gamma_1$).
  Let $h$ be a nearly Laplace Transform
$$ h(u)= \un{u>0}\, l(u)\, \intof e^{ -ux}\mu(dx)\,,$$
where $l$ is locally bounded and $\mu$ a complex measure such that
$$ \intof (1\wedge \unsur{\sqrt{x}})\valabs{\mu}(dx) < +\infty\,.$$
}
\begin{proposition}
  Suppose that $\mu$ is a sigma-finite measure on $(0,+\infty)$ such that for a $p>1$,
 $$\int \mu(dx)\,\sup(x^{-p/2}, x^{-1/2}) < +\infty\,.$$
 Then, for any measurable 
$\psi:L^1(\mu)\to \R$ such that $\esp{\valabs{\psi(\yf)}}<+\infty$
$$ \unsur{t} \intot \psi(\Yrond_s)\, ds \cvt \esp{\psi(\yf)}\,,$$
where the convergence is almost sure and in $L^1$.
In particular,  $h(u)=\intof e^{-ux}\,\mu(du)$ is in $L^2(0,+\infty)$; hence, if  we let $N$ denote a standard Gaussian random variable, $a=\norme{h}_{L^2(\rp)}=\etp{\intof h(u)^2\, du}^\undemi$, and 
$$ V_h(t) = \intot h(t-u)\, dB_u\,,$$
 then, for any measurable $\psi$ such that $\esp{\valabs{\psi(N)}}< +\infty$, for any $y\in L^1(\mu)$, under $\PP_y$,
$$ \unsur{t} \intot \phi(\unsur{a} V_h(s))\, ds \cvt \esp{\phi(N)}$$
where the convergence takes place almost surely and in $L^1$.
\end{proposition}
\begin{proof}
  The second part of the proposition is an obvious consequence of the first part: we only need to take $\psi(y)=\phi(\int y(x)\,d\mu(x))$ and to use Fubini--Tonelli's theorem
  \begin{align*}
    \intof h(u)^2\, du &= \intof du \etp{\int e^{-xu}\, d\mu(x)}\etp{\int e^{-yu}d\mu(y)} \\
&= \int\int d\mu(x)\,d\mu(y)\, \intof e^{-(x+y)u}\, du\\
&= \int \int d\mu(x)d\mu(y)\, \unsur{x+y}\,.
  \end{align*}

To prove the first part, we consider the invariant sigma-field $\Irond$,

$$ F=\ens{z:t \in \rp \to z_t \in L^1(\mu),\, \unsur{t} \intot \psi(z_s)\, ds \cvt \esp{\psi(\yf)\mid \Irond}}\,,$$
and $f(y)=\PP_y\etp{(t\to \Yrond_t)\notin F}$. Birkhoff's ergodic theorem implies that $\esp{f(\yf)}=0$ (and that the convergence also takes place in $L^1$).
The invariance of $f$ is an easy consequence of the Markov property of $(\Yrond_t\tgo)$ since  for all $t\ge 0$: $F= \theta_t^{-1}(F)$ where $\theta_t$ denotes the shift operator. Proposition~\ref{mar:pro:cvloiy} implies that $f$ is a constant, and thus, for all $y\in L^1(\mu)$, $f(y)=0$.

All we need to show now is that the invariant sigma-field $\Irond$ is almost surely trivial. Let $Z$ be a bounded $\Irond$ measurable random variable. We have, for all $t$, $Z=Z\circ \theta_t$ and therefore the function $f(y)=\E_y\etc{Z}$ is constant since

\begin{equation}
\label{rely}
 \esp{Z \mid \Frond_t}=\esp{Z\circ \theta_t \mid \Frond_t}= f(\Yrond_t)\,.
\end{equation}
Taking expectations under $\PP_y$ yields 
$$P_tf(y)=\E_y\etc{f(\Yrond_t)}=f(y)\,.$$
Proposition~\ref{mar:pro:cvloiy} implies that $f$ is constant, say $f(y)=c$ for every $y\in L^1(\mu)$. Injecting this in relation~\eqref{rely} yields, that under $\PP_y$
$$ \esp{Z \mid \Frond_t}=c \qquad(t\ge 0)\,.$$
Taking limits as $t\to +\infty$ gives finally that $Z=c$ a.s.

\end{proof}

\providecommand{\bysame}{\leavevmode\hbox to3em{\hrulefill}\thinspace}

\end{document}